\documentclass[12pt]{article}
\usepackage[latin2]{inputenc}
\usepackage{amssymb}
\usepackage{amsmath,amsopn}
\usepackage{amsthm}
\usepackage{authblk}
\usepackage{bm}
\usepackage{color}
\usepackage{enumerate}
\usepackage{epsfig,latexsym}
\usepackage{epstopdf}
\usepackage{graphicx}
\usepackage{multirow}
\usepackage[authoryear,sort,round]{natbib}
\usepackage{times}
\usepackage{t1enc}
\usepackage{url}

\newcommand{\varr}{\! V\! \! A\! R}
\begin{document}

\setcounter{page}{1}
\newcommand\balline{\small László Varga and András Zempléni}
\newcommand\jobbline{\small Generalised block bootstrap and its use in meteorology}

\title{Generalised block bootstrap and its use in meteorology}

\author[1]{László Varga\thanks{Corresponding author}\thanks{vargal4@math.elte.hu}}
\author[1]{András Zempléni} %vl0420 én benne hagynám az E-mail címedet \thanks{zempleni@math.elte.hu}}
\affil[1]{Department of Probability Theory and Statistics, E\"otv\"os Lor\'and University, Budapest, Hungary}

\maketitle

\abstract { 
In an earlier paper \cite{rvz}, we have emphasized the effective sample size for autocorrelated data. The simulations were based on the block bootstrap methodology. However, the discreteness of the usual block size did not allow for exact calculations. In this paper we propose a generalisation of the block bootstrap methodology, relate it to the existing optimisation procedures and apply it to a temperature data set.
Our other focus is on statistical tests, where quite often the actual sample size plays an important role, even in case of relatively large samples. This is especially the case for copulas. These are used for investigating the dependencies among data sets. 
As in quite a few real applications the time dependence cannot be neglected, we investigated the effect of this phenomenon to the used 
test statistic. %vl0420 test-statistics.
The critical values can be computed by the proposed new block bootstrap simulation, where the block sizes are  
determined e.g. by fitting a VAR model to the observations.
The results are illustrated for models of the used temperature data. \\[1mm]
\begin{tabular}{lp{90mm}}
\hspace{-2.5mm}\textbf{Keywords:} & 
block bootstrap, copula, temperature data, test of homogeneity, VAR models
\end{tabular} 
}  % abstract vége

%%% 1. fejezet %%%

\section{Introduction}
In the last decades the bootstrap methodology has become  more and more widespread in different areas of statistical applications. See e.g. \citet{ch} for a review of possible areas from spatial models to financial data and data mining, where bootstrap may be used. In this paper we focus on the effect of the serial dependence, naturally arising in many time series data. The bootstrap samples must match the dependence within the data, so the block bootstrap is the suggested method for bootstrapping time series. \citet{hh} investigate this approach in some detail, including suggestions for selecting the optimal block size.
In an earlier paper \citet{rvz}, we investigated the possibilities for using the block bootstrap methods for checking the validity of the copula models. In this paper we present an improvement to the classical block bootstrap methodology, which is especially relevant in our applications. 

In Section 2, we first briefly review the elements of stationary time series. 
In the bivariate case, the vector autoregression (VAR) process is one of the most important models, becoming popular first in the area of econometrics (\citet{si}).
%originally developed in the area of econometrics ().  
For recent applications of VAR models in meteorology, see for example \citet{hbm}, \citet{nyl} or \citet{fa}.
We briefly present the main properties of VAR models, which are used in the sequel and present the notations. 

In Section 3, we introduce the concept of copulas, 
the most convenient objects for analysing the dependence structures among variables. Their history goes as back as \citet{ho}, but their applications are much more recent. However, they have spread very quickly to the most important areas -- for a recent analysis in meteorology, see \citet{cb}. Most of these works 
use different parametric copula models, but we are more interested in testing for possible changes in the dependency structure of our temperature data,
%vl0414 including some of the most important copula families  
so also introduce the most recent approaches in testing homogeneity of such models, which are based on the empirical copula process.
% -- testing the goodness of fit (GoF) for such models. 

Section 4 is devoted to the bootstrap resampling method, including the block bootstrap approach, which is suitable for the case of serial dependent observations. Here we introduce a generalisation, which helps overcoming the problem that 
originally the block size was supposed to be a natural number. In our approach the block size is a random variable, but it 
%has such a small range that it 
contains the original block bootstrap as a special case. Due to this small variance in the sample size, it overcomes the problem of extensive random error in case of the so-called stationary bootstrap (\citet{pr}).

In Section 5 we apply our approach to the gridded temperature data base of E-OBS, which is a product of the EU-FP6 project ENSEMBLES (\citet{hhgk}). Here we use the daily mean temperature data from the 0.5 grade-grid. 
%We check  the effect of the serial dependence on the model selection and compare the possible approaches to block-size determination. 
%za0419
Our focus is on checking for possible changes in the dependence pattern between the grid point of Budapest and some other grid points within the Carpathian Basin. We show that in some
cases there is a significant deviation from homogeneity of the first and second part of the data.
%za0419 pontosítani kellene
The conclusion summarizes our findings and gives some interesting open questions. 
    % bevezetes

%%% 2. fejezet %%%
%time series

\section{Vector autoregression (VAR) processes}
\label{secVAR} 

We call the $d$-dimensional series $\{\mathbf{X}_t \} _{t \in \mathbb{Z}} = (\mathbf{X}_0, \mathbf{X}_{\pm 1},\mathbf{X}_{\pm 2},\ldots)$ 
a time series if its elements are $d$-dimensional random vectors, which are usually not independent from each other.
%za0525,
Here we consider $t$ as the time. 
Let us assume that the random variables have finite second moments. The time series $\{ \mathbf{X}_t \} _{t \in \mathbb{Z}}$ is weakly stationary 
%za0525
(later we just call it stationary), if for all $t$, if neither the mean function $E(\mathbf{X}_t)$ 
%is a constant vector and 
nor the covariance matrix Cov$(\mathbf{X}_{t+s},\mathbf{X}_t)$ 
%does not 
depends on $t$ for all $s\in \mathbb{Z}$. Stationarity is an important property, it means that the time series is translation invariant. 

One of the most frequently applied time series models are the so-called autoregressive (AR) processes and 
%za0515 its 
their multidimensional counterparts, the vector-autoregressive (VAR) models. In the following, we define the VAR($p$) process and give its main properties in two dimensions as this is necessary 
%za0525for 
to our applications. 

The time series $\{ \mathbf{X}_t \} _{t \in \mathbb{Z}}=\{ (X_{1,t},X_{2,t})^T \} _{t \in \mathbb{Z}}$ is called 
%za0515
a zero-mean two-dimensional VAR($p$) process if
\begin{equation}
\label{var}
 \mathbf{X}_t=A_1 \mathbf{X}_{t-1}+A_2 \mathbf{X}_{t-2}+\ldots+ A_p \mathbf{X}_{t-p}+\bm{\varepsilon}_t,
\end{equation}
where $A_1,\ldots,A_p$ are $2\times 2$ parameter matrices and the  $\{ \bm{\varepsilon}_t  \} _{t \in \mathbb{Z}}$ independent innovation process is a two-dimensional white noise with $E(\bm{\varepsilon}_t)=\mathbf{0}=(0,0)^T$ and Cov$(\bm{\varepsilon}_t)=C$
%vl0517 :=\text{diag}(\sigma_1^2,\sigma_2^2)$ diagonal 
symmetric positive definite covariance matrix. The VAR($p$) process is stationary if the roots of the $P(x)=\det(I_2-A_1x-\ldots-A_px^p)$ characteristic polynomial lie outside the unit circle. 

Any VAR($p$) process can be rewritten as a VAR(1) process in the following way: $\mathbf{Y}_t=A \mathbf{Y}_{t-1}+\mathbf{e}_t$, where \\
$\mathbf{Y}_t=\begin{pmatrix}
	\mathbf{X}_t \\ \mathbf{X}_{t-1} \\ \vdots \\ \mathbf{X}_{t-p+1}
\end{pmatrix}$,
$\mathbf{e}_t=\begin{pmatrix}
\bm{\varepsilon}_t \\ \mathbf{0} \\ \vdots \\ \mathbf{0}
\end{pmatrix}$ and 
$A=\begin{pmatrix}
A_1        & A_2        & \cdots & A_{p-1}    & A_p \\ 
I_2        & \mathbf{0} & \cdots & \mathbf{0} & \mathbf{0} \\
\mathbf{0} & I_2        & \cdots & \mathbf{0} & \mathbf{0} \\ \vdots     & \vdots     & \ddots & \vdots     & \vdots \\
\mathbf{0} & \mathbf{0} & \ldots & I_2        & \mathbf{0}
\end{pmatrix}$.
This representation is more convenient in calculating the autocovariances.
An equivalent condition for stationarity is that all the eigenvalues of the coefficent matrix $A$ are smaller than one in modulus. In this case the time series is causal -- it has an infinite moving average representation in the form $\mathbf{Y}_t= \sum\limits_{i=0}^{\infty} A^i \mathbf{e}_{t-i}  $.
%vl0517 , which is really helpful in calculating the autocovariances

For the remainder of this section we assume that $\mathbf{X}_t$ is stationary. Let us denote with $\Gamma_X(h)=E(\mathbf{X}_{1+h} \mathbf{X}_1^T)$ the autocovariance function of the process $\mathbf{X}_t$. $\Gamma_X(h)$ is a $2\times 2$ matrix valued function, the symbols $\gamma_{i,j}(h)$ stand for its elements. We denote with $\Gamma_Y(h)=E(\mathbf{Y}_{1+h} \mathbf{Y}_1^T)$ the $2p\times 2p$ matrix valued autocovariance function of the process $\mathbf{Y}_t$. The covariance matrix of $\mathbf{Y}_t$ is $\Gamma_Y(0)$
which can be determined by solving the matrix equation $\Gamma_Y(0)-A\, \Gamma_Y(0) A^T= \text{Cov} (\mathbf{e_t})$. It is easy to see that 
for $1 \leq h \in \mathbb{Z}$, the autocovariances can be calculated by $\Gamma_Y(h)=A^{h}\Gamma_Y(0)$. The powers of the matrix $A$ 
%are easily derived with 
can easily be computed using 
the spectral decomposition. 
Lastly, we need the autocovariance matrix of the original process, and by the construction, it is the upper left $2 \times 2$ submatrix of $\Gamma_Y(h)$.

In the applications we will use the covariance matrix of the sample mean. The following asymptotic result will be crucial in our investigations: if $\sum\limits_{h=-\infty}^{\infty} |\gamma_{i,i}(h)|<\infty $ for $i=1,2$, then
\begin{flalign}
\label{elmtr}
n \cdot \text{tr} \left( \text{Cov}(\overline{\mathbf{X}}_n) \right) \underset{ }{\longrightarrow} \sum\limits_{i=1}^2\sum\limits_{h=-\infty}^{\infty} \gamma_{i,i}(h) \ \ \ \ \text{as } n \to \infty,
\end{flalign}
where tr$(\cdot)$ denotes the trace of a matrix.

It is important to check if the chosen time series model is adequate. If the model fits well, the fitted residuals should behave as a realisation of a white noise process. The hard part is to check whether the residuals are independent, thus there is no serial dependence. There are several methods for verifying this problem, the most standard is the Ljung-Box test, which tests whether a specified group (usually the first 10-20 lags) of autocorrelations is different from zero. Another often applied serial correlation test is the Breusch-Godfrey test. A more recent multidimensional approach was published in the paper \citet{ky}, where the test was based on the empirical copula process. The main ideas and the concept of this test stems from \citet{gr}.

For further details about time series see for example \citet{bd} or \citet{ss}.

    % idősorok, AR folyamatok
%%% 3. fejezet %%%
%%GoF and copula

\section{Copulas and their goodness-of-fit}
\label{ChCop}

Let $\mathbf{X}=(X_1,\ldots,X_d)^T$ be a random vector with joint distribution function %${H(x_1,\ldots,x_d)=F_{X_1,\ldots,X_d}(x_1,\ldots,x_d)}$ 
$F_{\mathbf{X}}(\mathbf{x}) = F_{X_1,\ldots,X_d}(x_1,\ldots,x_d)$ and marginal distribution functions \mbox{$F_1(x_1)=F_{X_1}(x_1),$ $\ldots, F_d(x_d)=F_{X_d}(x_d)$}. 
Sklar's theorem claims that there exists a copula ${\bf C}$, a distribution over the $d$-dimensional unit cube, with uniform margins, such that
\begin{equation*}
%H(\mathbf{x})={ H}(x_1,\ldots,x_d)=
F_{X_1,\ldots,X_d}(x_1,\ldots,x_d)=
{\bf C}\left(F_1(x_1),\ldots,F_d(x_d)\right).
\label{joint}
\end{equation*}
Moreover the copula $\mathbf{C}$ is unique if the marginal distribution functions are continuous. This construction allows the investigation of the dependence structure without specifying the marginal distributions. In the recent literature various families of copulas have been introduced, for an overview and examples see e.g. the introductory textbook of \cite{ne}.

In this paper we focus on testing the homogeneity of copulas, motivated by the question whether the climate change has also an effect on the dependence between pairs of observations. If this change is indeed observable, then it may have a substantial effect on the spatial structure of temperature anomalies, worth for further meteorological investigations. So we do not have to go into the parametric inference, as we are just interested in the homogeneity analysis.

%\subsection{Homogeneity tests for copulas} 
%\textbf{Homogeneity tests for copulas}

\begin{figure}[Ht!]
	\centering{
		\includegraphics[width=60mm]{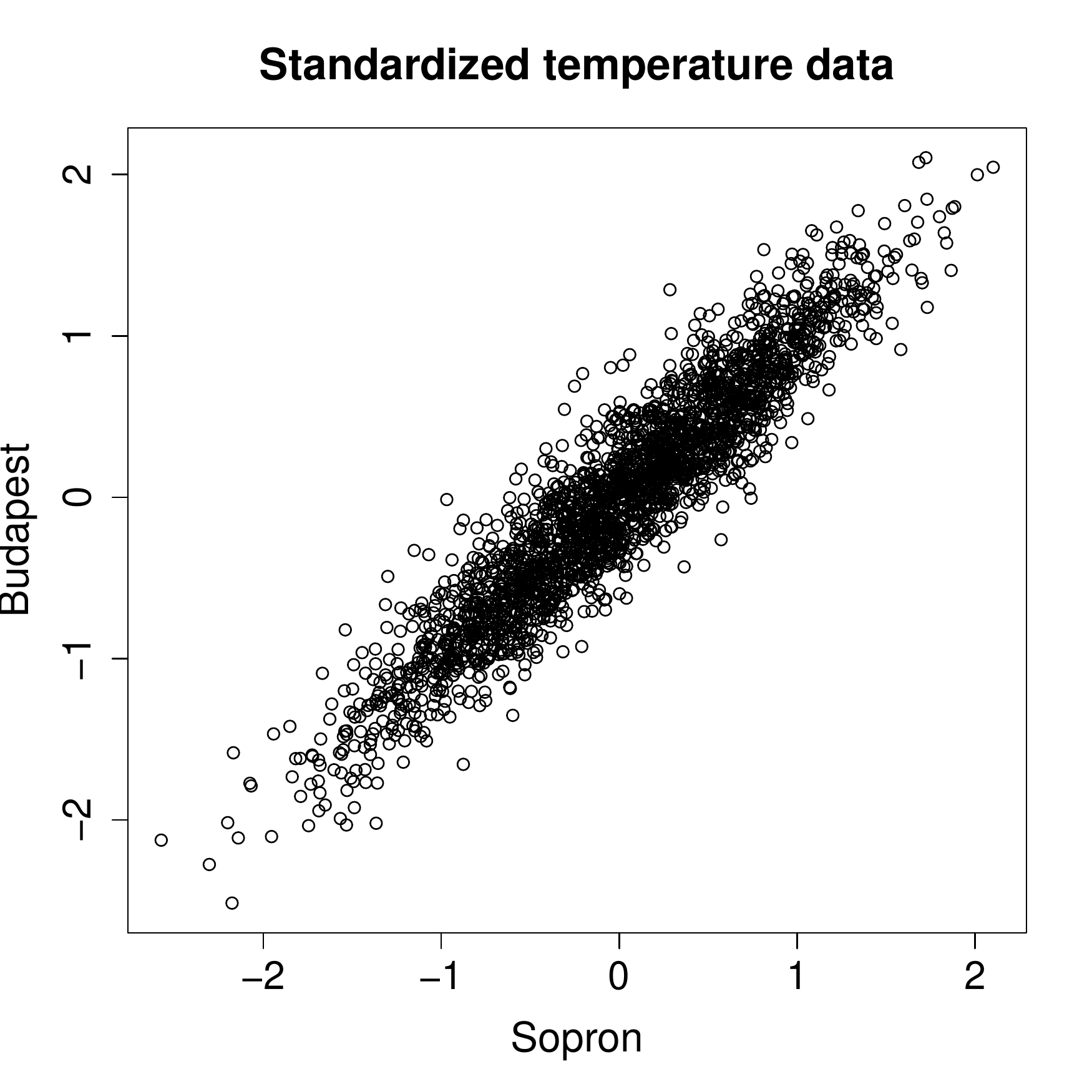}
		\includegraphics[width=60mm]{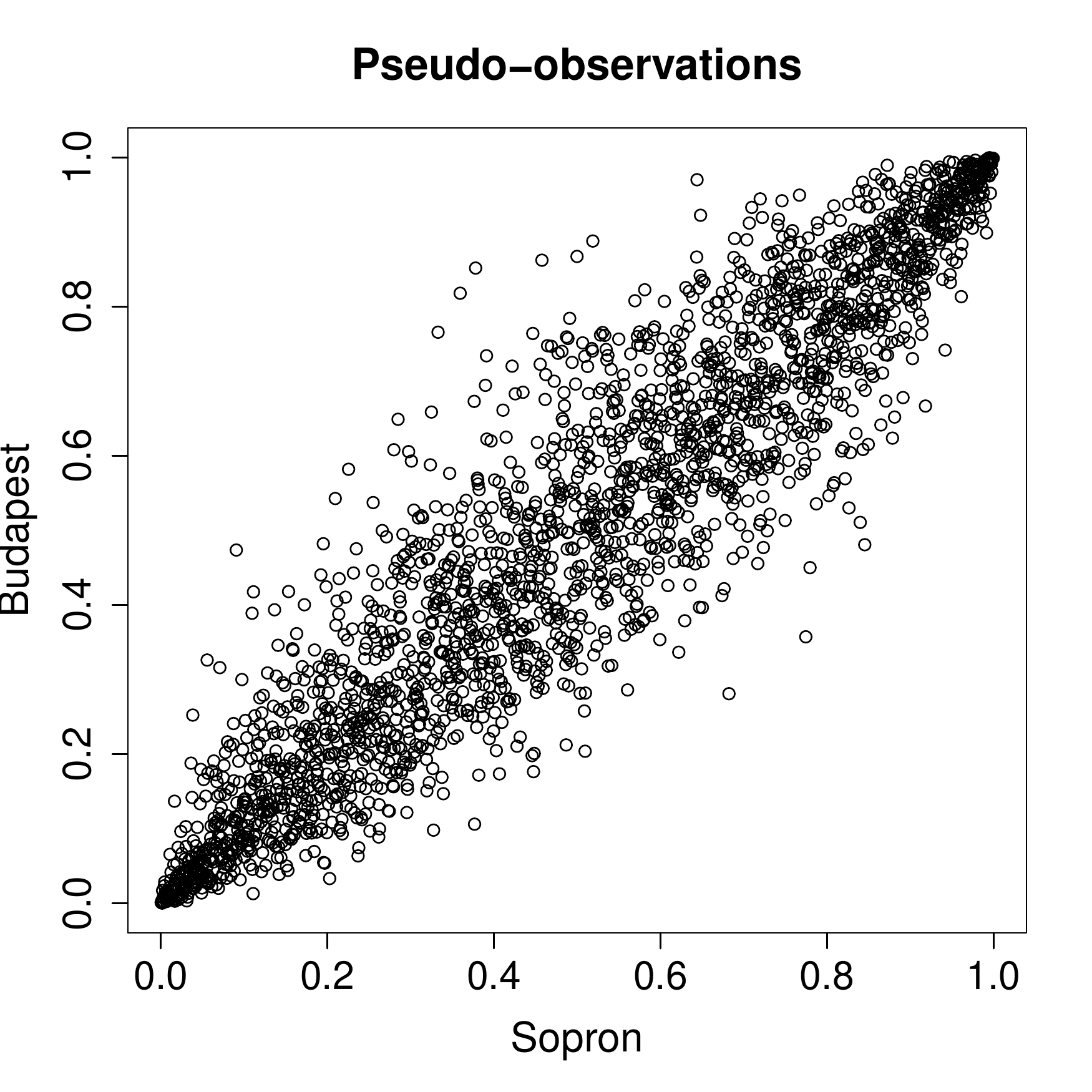}
	}
	\caption{Bivariate data and the corresponding pseudo-observations }
	\label{kfun}
\end{figure}

Let us suppose we have two independent samples of $\mathbb{R}^d$-valued vectors. The first sample is $\mathbf{X}_1,\ldots,\mathbf{X}_n$ and the second one is $\mathbf{Y}_1,\ldots,\mathbf{Y}_m$.
Formally we intend to test the hypothesis that the dependence structure of the two copulas has arisen from the same copula $\mathbf{C}_0$. The most obvious way for 
%testing the GoF between two copulas 
testing the homogeneity of two copulas
is to consider multidimensional $\chi^2$ approaches, but in this case we need to discretize the data, losing valuable information. 
In order to avoid its use, we can follow the approach of \cite{rs}, who have developed a method for this problem. Their approach is based on the empirical copula, defined for the first sample as 
\begin{equation*}
\label{kt}
\mathbf{C}_{1,n}(\mathbf{u})=\frac{1}{n}\sum_{i=1}^n I(\mathbf{U}_i\le \mathbf{u}). 
\end{equation*}
where $\mathbf{u}\in \mathbb{R}^d$ and $\mathbf{U}_i$ denotes the $d$-dimensional vector of the rank based pseudo-observations: $\mathbf{U}_i=\mathbf{U}_{i,n}=(U_{i1,n},\ldots,U_{id,n})$, where $n$ refers to the size of first sample and $U_{ij,n}=\frac{n}{n+1}F_{j}(X_{ij}) $.
For illustrations see \hbox{Figure \ref{kfun}}. 
Similarly, based on the pseudo-observations $\mathbf{V}_i$ of the second sample, we can define the empirical copula  $\mathbf{C}_{2,m}(\mathbf{u})$.
%We define$\mathbf{C}_{2,m}(\mathbf{u})$ %similarly with $\mathbf{V}_i$ representing the vector of the rank based pseudo-observations for the second sample.

The proposed tests for checking the homogeneity of two samples are based on functionals of the
empirical process:
\begin{equation*}
\kappa_{n,m}(t)=\frac{\mathbf{C}_{1,n}(t)-\mathbf{C}_{2,m}(t)}{\sqrt{\frac{1}{n}+\frac{1}{m}}},
\end{equation*} where the asymptotic
properties of the statistic can be based on the limit of the empirical copula processes. 
There are two different kind of approaches investigated in \citet{gen}: 
%--  in case of goodness-of-fit --
the Cram\'er-von Mises
type statistic $S_{n,m}=\int\limits_{0}^{1}\left( \kappa_{n,m}(t) \right) ^{2}dt$, and the Kolmogorov-Smirnov type statistic
$T_{n,m} = \underset{0 \leq t \leq 1}{\sup} \left| \kappa_{n,m}(t) \right|$. As the second approach is considered
to be
generally less powerful, we based our inference on the  statistic $ K_{*}=\frac{1}{N}\sum\limits_{t_{i}\in (0,1)}\left( \kappa_{n,m}(t_i) \right)^2$,
where $(t_{i})_{i=1}^N$ is an appropriately fine division of the interval $(0,1)$. 
After some %exposition, 
calculations,
the Cramér-von Mises test statistic can be written in the following form (see \cite{rs}):
\begin{flalign}
\label{tstat}
S_{n,m}=\left( \frac{1}{n}+\frac{1}{m} \right)^{-1} \cdot &
\left[ 
\frac{1}{n^2} \sum\limits_{i=1}^n \sum\limits_{j=1}^m \prod\limits_{s=1}^d (1-U_{is,n} \vee U_{js,n}) +  \right.  \\ 
& \hspace{3mm} + \frac{1}{m^2} \sum\limits_{i=1}^m \sum\limits_{j=1}^m \prod\limits_{s=1}^d (1-V_{is,m} \vee V_{js,m} - \nonumber \\ 
& \hspace{3mm} - \left. \frac{2}{nm}\sum\limits_{i=1}^n \sum\limits_{j=1}^m \prod\limits_{s=1}^d (1-(U_{is,n} \vee V_{js,m}) 
\right], \nonumber 
\end{flalign}
where $u \vee v =\max(u,v)$.

As the limit distribution of the above statistic is not distribution-free, a simulation algorithm is needed to get critical values. The algorithm is the following:
\begin{enumerate}
\item Generate an $n$-element bootstrap sample from the first of the observation vectors -- see its details in the next section --, and compute the statistic $S_{n,m}$, based on this and the original first sample; 
\item Repeat the above steps as many times as needed to an accurate estimation of the $p$-values. The distribution function $G$ of the bootstrap statistic can be considered as a good approximation of the test statistic under $H_0$, so we used $1-G(S_{n,m})$ as the estimate of the $p$-value.
%, as now we reject the large values of our statistic.
\end{enumerate}

    % copulák, homogenitásvizsgálat

%%% 4. fejezet %%%
\section{Bootstrap methods}
The bootstrap is a usually computer-intensive, resampling method for estimating the distribution of a statistic of interest.
The concept of the bootstrap was introduced in the classical article by Ben Efron (\cite{ef}) and since then, it has become one of the most widely used Monte Carlo methods in a number of aeras of applied sciences.

\subsection{Bootstrap for i.i.d. data}

Let $\mathcal{X}_n=(X_1,\ldots, X_n)^T$  be a sequence of independent, identically distributed (i.i.d.) random variables with unknown common univariate distribution $F$ and let $T_n=t_n(\mathcal{X}_n; F)$ be a statistic (like the sample mean $\overline X$). Our main purpose is to approximate the unknown distribution of $T_n$ or its function of interest, for example its standard deviation (the standard error) or some quantiles. 

The basic bootstrap method (mostly referred as i.i.d. bootstrap) is the following. For a given $\mathcal{X}_n$, we
draw a random sample {$\mathcal{X}_m^*=\{X_1^*,\ldots,X_m^*\}$ of size $m$ (usually $m = n$) with replacement from $\mathcal{X}_n$. Therefore, the common distribution of the $X_i^*$'s is given by the empirical distribution $\hat{F}_n=n^{-1}\sum_{i=1}^n\delta_{X_i}$, 
where $\delta_z$ is the probability measure having unit mass at $z$. In the next step, we define the bootstrap version of the statistic $T_n$: \mbox{$T_{m,n}^*=t_m(\mathcal{X}_m^*;\hat{F}_n)$}. By repeating this procedure, we can approximate the unknown distribution $G_n$ by its bootstrap counterpart $G_n^*$. In most of the cases the distribution of $G_n^*$ cannot be determined explicitely, but it can be approximated by simulation. 
	
\subsection{Block bootstrap methods}
	
In our case we are interested in the effect of serial dependence on the homogeneity tests and on modelling in general, for example 
on the covariance matrix of our estimators. 
If the data are dependent then the estimates based on i.i.d. bootstrap methods my not be consistent.

In the presence of serial dependence, one of the most commonly used methods is the so-called block bootstrap, see \cite{lah} for details. In this paper, we generalise the circular block
bootstrap (CBB) which can be defined as follows. First, we wrap the data $X_1,\ldots,X_n$ around a circle, i.e., define the series
$Y_t=X_{t_{\rm{mod}(n)}}$ 
($t\in\mathbb{N}$), where $\rm{mod}(n)$ 
denotes division "modulo $n$". 
This means that $X_k=Y_{k}=Y_{k+n}=Y_{2k+n}=\ldots$ for all $k\in \{ 1,2,\ldots,n\}$. %vl0410 
For some $m$, let $i_1,\ldots,i_m$ be a uniform sample from the set $\{1,2,\ldots,n\}$.
Then, for a given block size $b$, we construct $n'=m\cdot b$ \ \ ($n'\approx n$) pseudo-data:
	
\begin{center} 
	$Y^*_{(k-1)b+j}=Y_{i_k+j-1}$  \qquad  where $j=1,\ldots,b$ \ \ and $k=1,\ldots,m$. 
\end{center}
Finally, let us calculate the function of interest, for example the bootstrap sample mean as follows:
$\overline{Y}^*_{n'}=\frac{Y^*_1+\ldots+Y^*_{n'}}{n'}$. 
	
Block length plays an important role in the process, and it is not trivial to determine its optimal value. 
Politis and White in their article \cite{pw} suggest
an "automatic" block length selection algorithm (its correction was published in \cite{ppw}) - but the practical applications of this method are far from obvious as there are parameters in it, which have to be chosen. 
	
We used a similar approach in our paper \cite{rvz}.
Our idea was that we tried to find the best block size by fitting a model and then checking the variance of $\overline X$ with the help of the block bootstrap. 
The block size was determined as the $\widehat{b}$, for which the estimated trace of the covariance matrix was the 
nearest 
to the one got by the fitted VAR model:
\begin{equation} \widehat{b}=\underset{1\leq b \,\in \mathbb{Z}}{\text{argmin}}
\left| 
\text{tr}\left( \text{Cov}\left( \overline{X}_{\varr} \right) \right) - 
\text{tr}\left( \text{Cov}_* (\overline{X}_b^* ) \right)
\right| ,\label{bopt}
\end{equation}
where $ \text{Cov}_* (\overline{X}_b^*)= \text{Cov} (\overline{X}_b^* | \mathcal{X}_n)$. In the literature, simulations are naturally based on integer block sizes. But using the block length of (\ref{bopt}), 
the estimated trace of covariance may be not be close enough to the theoretical trace of covariance.
The same is true for other methods for block size determination.
This may cause substantial bias, as in our case the relative difference between subsequent values of tr$(\text{Cov}_* (\overline{X}_b^* ))$ can be quite large, especially for small $b$. This can be overcome by the following generalisation of the block bootstrap methodology. 
	
In case of $b\!>\!1$, \ $b\!\in \! \mathbb{R}$, let the generalised block bootstrap sample be defined as follows. Let $k$ be a random integer between $1$ and the sample size $n$, and again, let us wrap the sample around the circle.
The bootstrap blocks are either of length $\lfloor b \rfloor$ or $\lceil b \rceil$: \\
\begin{tabular}{llr} 
	$X_k,X_{k+1},\dots,X_{k+\lfloor b \rfloor}$ & with probability & $1- b + \lfloor b \rfloor $  \\ 
	$X_k,X_{k+1},\dots,X_{k+\lceil b \rceil}$ & with probability &  $b- \lfloor b \rfloor$
	% \label{bdef}
\end{tabular} \\
where $\lceil b \rceil$ denotes the upper, and $\lfloor b \rfloor$ the lower integer part of $b$. 
At last, we put the blocks together. This procedure ensures that for integer-valued $b$ the new definition coincides with the traditional one, so this is indeed a generalisation. In the applications (Section \ref{secapp}) we show the clear advantages of this approach. 
Actually all of the relevant algorithms for finding the optimal block size can easily be adapted to find a solution in this generalised sense. 
In our case, instead of (\ref{bopt}), we simply solve the equation
\begin{equation}
 \text{tr}\left( \text{Cov} ( \overline{X}_{\varr} ) \right) = \text{tr} \left( \text{Cov}_* (\overline{X}_b^* ) \right) .\label{bopt2}
\end{equation}   

In the same way as the circular block bootstrap sample, our generalized bootstrap sample is not a stationary process, conditional on the original sample. The block bootstrap sample is in only one case conditionally stationary: when the block lengths follow a geometric distribution, independent from each other (\citet{pr}).

The covariance matrix $ \text{Cov}_* (\overline{X}_b^*)$ can be explicitly calculated. Let us denote with $L_1,L_2,\ldots$ the block sizes -- they are random variables independent from each other with common conditional distribution\footnote{$P_*$ and $E_*$ denotes the conditional distribution and the conditional expectation given the sample $\mathcal{X}_n$, so $P_*(L_1= \lceil b \rceil)=P(L_1= \lceil b \rceil | \mathcal{X}_n)$ and $E_*(L_1)=E(L_1 | \mathcal{X}_n)$. }
 $P_*(L_1=\lceil b \rceil) = 1-P_*(L_1=\lfloor b \rfloor)=b-\lfloor b \rfloor$. We can also write $L_i=\lfloor b \rfloor+J_i$, where $J_i| \mathcal{X}_n$ follows a Bernoulli distribution with parameter $p=b-\lfloor b \rfloor$. 
Let $N$ be the random variable, which gives the number of blocks with block size $\lfloor b \rfloor$. If we have $N$, we can calculate the number of blocks with block size $\lceil b \rceil$, we denote it with $g(N)$. So $g(N)=\left\lfloor \frac{n-N\cdot \lfloor b \rfloor}{\lceil b \rceil} \right\rfloor$. Let us denote the remainder block size with $r(N)$, we can calculate it from the others: $r(N)=n-N\cdot \lfloor b \rfloor-g(N)\cdot \lceil b \rceil$. It can be seen that the conditional covariance matrix of the bootstrap mean can be calculated the following way:
\begin{flalign} 
\label{blvar}
\text{Cov}_* (&\overline{X}_b^*) = \frac{\lfloor b \rfloor^2}{n^2} \left[ \text{Cov}_*(\overline{X}_{\lfloor b \rfloor,i}^*) \cdot E_*N+ \text{Cov}_*N \cdot \overline{X}_n (\overline{X}_n)^T \right] +  \\
& + \frac{\lceil b \rceil^2}{n^2} \left[ \text{Cov}_*(\overline{X}_{\lceil b \rceil,i}^*) \cdot E_*(g(N))+ \text{Cov}_*(g(N)) \cdot \overline{X}_n (\overline{X}_n)^T \right] + \nonumber \\
& + \frac{1}{n^2}  \left[ \sum\limits_{i=0}^{\lfloor b \rfloor} i^2 P_*(r(N)=i) \cdot \text{Cov}_* ( \overline{X}_{i,1}^*) + \text{Cov}_*(r(N)) \cdot \overline{X}_n (\overline{X}_n)^T  \right]. 
\nonumber 
\end{flalign}
At last, we have to mention that the Politis and White algorithm actually gives a real and not an integer as the optimal block size -- this could be used without any rounding by our proposed method. The main problem is that the algorithm gives an extremely large block length, making meteorologically no sense.
	
    % bootstrap

%%% 5. fejezet %%%
% alkalmazások

\section{Applications}
\label{secapp}

The used observations are the 63 years of daily temperature data of the European Climate Assessment (E-OBS, \url{http://www.ecad.eu}). The methodology of deriving the data for the grid points has been published in \citet{hhgk}, where this  database has been used extensively for climate analysis.
We have worked 
with the part of the 0.5-grade grid -- available for whole Europe and northern Africa -- which lies in the Carpathian Basin.
%with the 0.5-grade grid in the Carpathian Basin, available for whole Europe and northern Africa. 
Figure \ref{map} depicts the used grid points. For later reference, we chose the grid point next to Budapest, one grid point in the neightborhood of the Hungarian capital and four further grid points lying far from Budapest, in different directions.

The quality of the data has been evaluated e.g. in \citet{hhnj}, 
 and it turned out to be reliable for most of Central Europe. 
As we have used the grid points, belonging to the Carpathian Basin, this validates our results. 

\begin{figure}[!h] 
	\centering{
		\includegraphics[width=12cm]{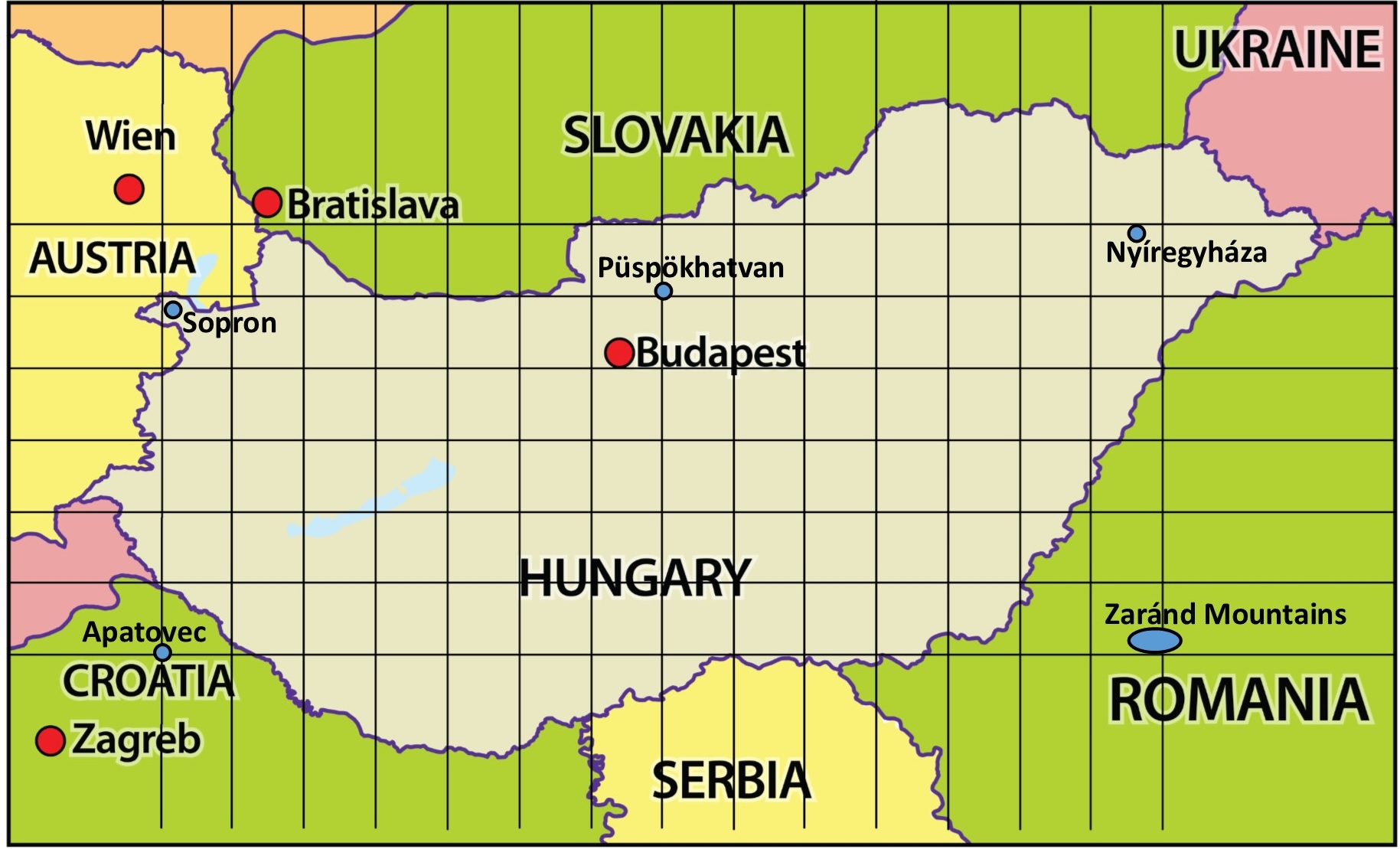}  	
	}
	\caption{The map of the Carpathian Basin with the used grid points.}
	\label{map}
\end{figure}

As we intend to use models, suitable for stationary data, first the stationarity had to be ensured. We have first subtracted the smoothed daily averages from the observations. The smoothing was made by loess regression, Figure \ref{daily}a.) and Figure \ref{daily}b.) depict the daily averages and standard deviations of the 63 years' data and the smoothing regression line for grid point Budapest.
It turned out that the second-order stationarity is still far from being true (in winter the variances were substantially larger than in summer), so we have divided the observations with the smoothed estimated standard deviation for the given day:
$$\tilde{x}_{t,n}=\frac{x_{t,n}-m_t}{s_t}$$
where $\tilde{x}_{t,n}$ is the standardized value for day $t$ in year $n$, based on the original observation $x_{t,n}$ for the same day and the smoothed average $m_t$ and smoothed standard deviation $s_t$.  Figure \ref{daily}c.) shows the original daily observations and the standardized data between January 1, 2010  and December 31, 2012 for grid point Budapest.

\begin{figure}[!h]
	\centering{
		    \includegraphics[width=6cm]{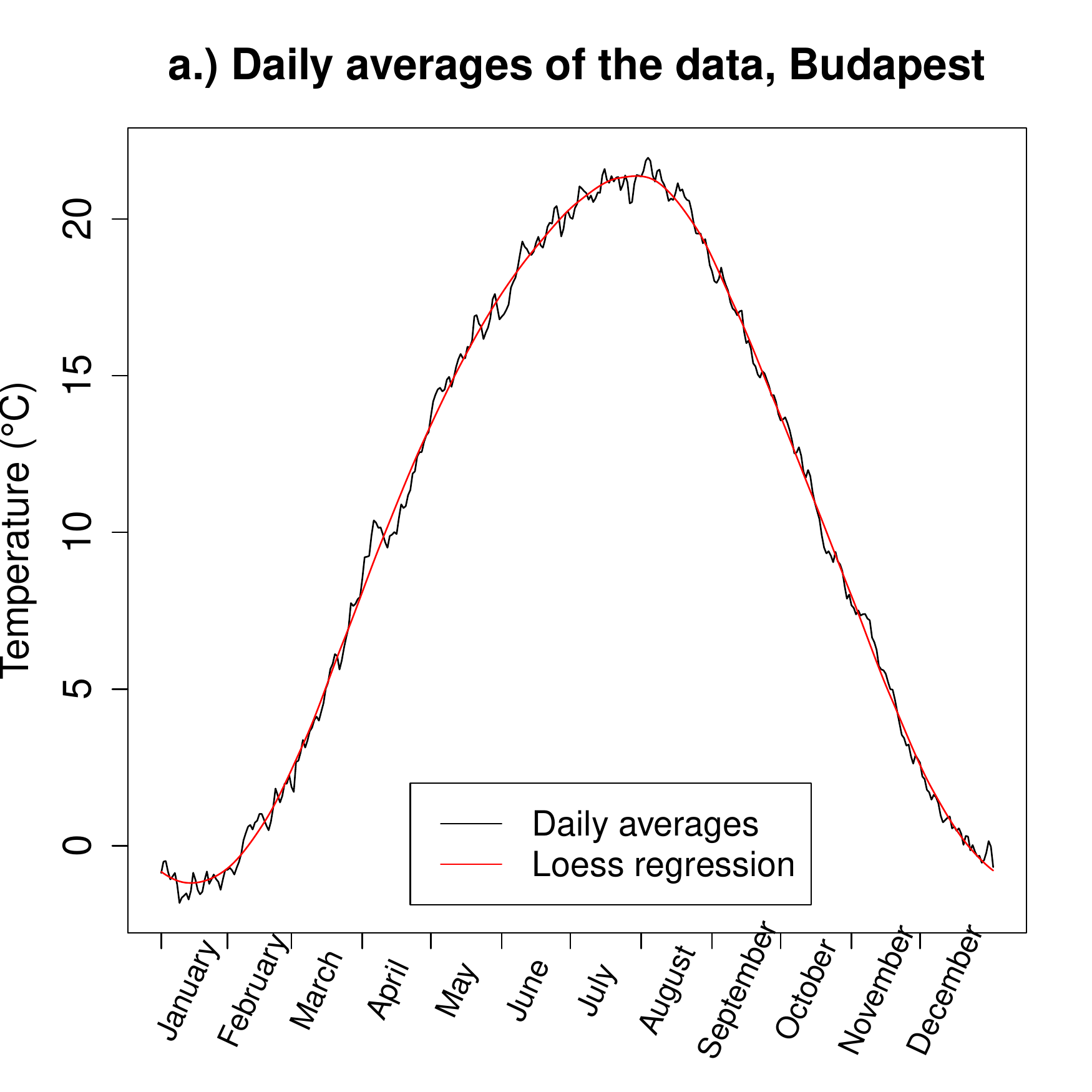} 
		    \includegraphics[width=6cm]{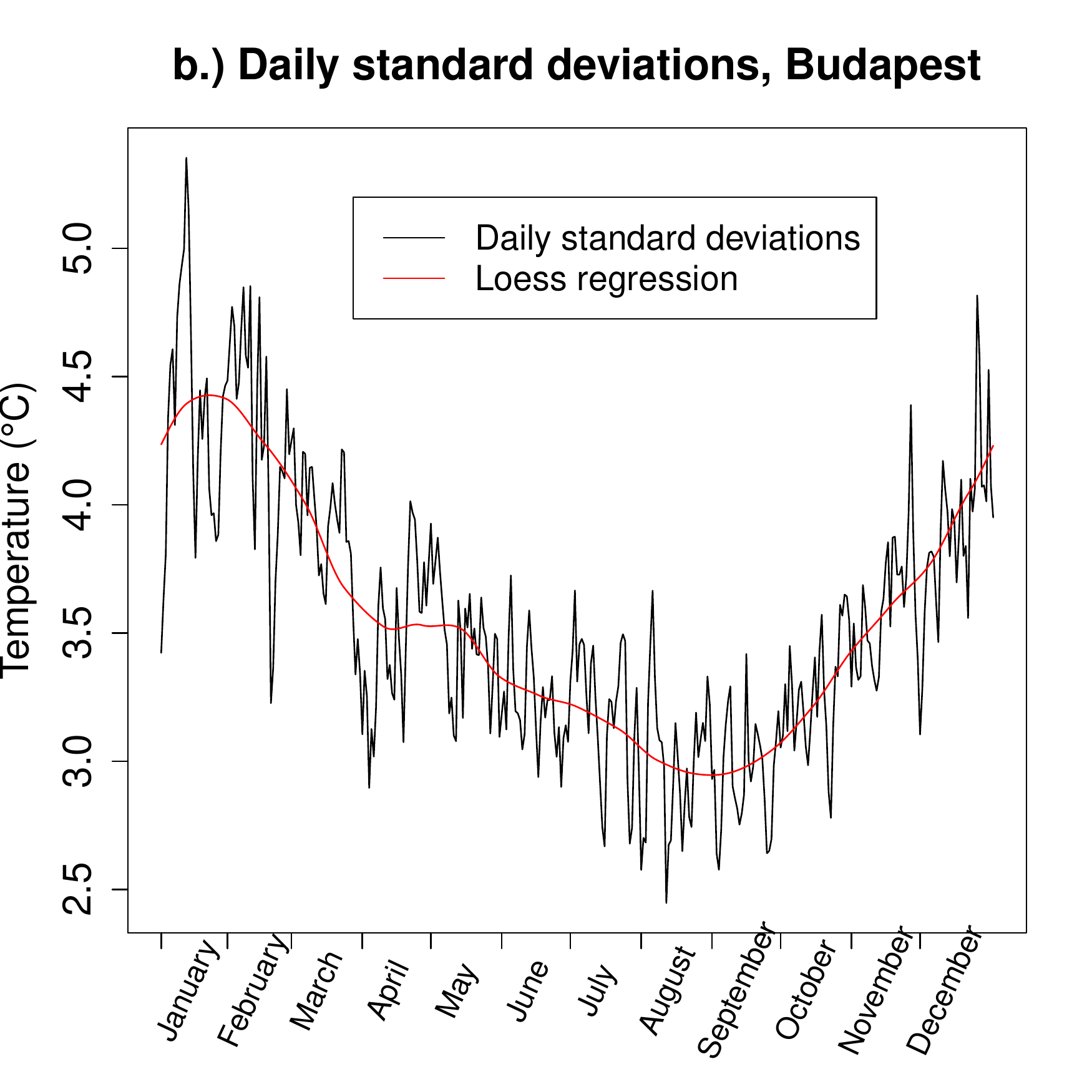} \\
    	     \includegraphics[height=6cm]{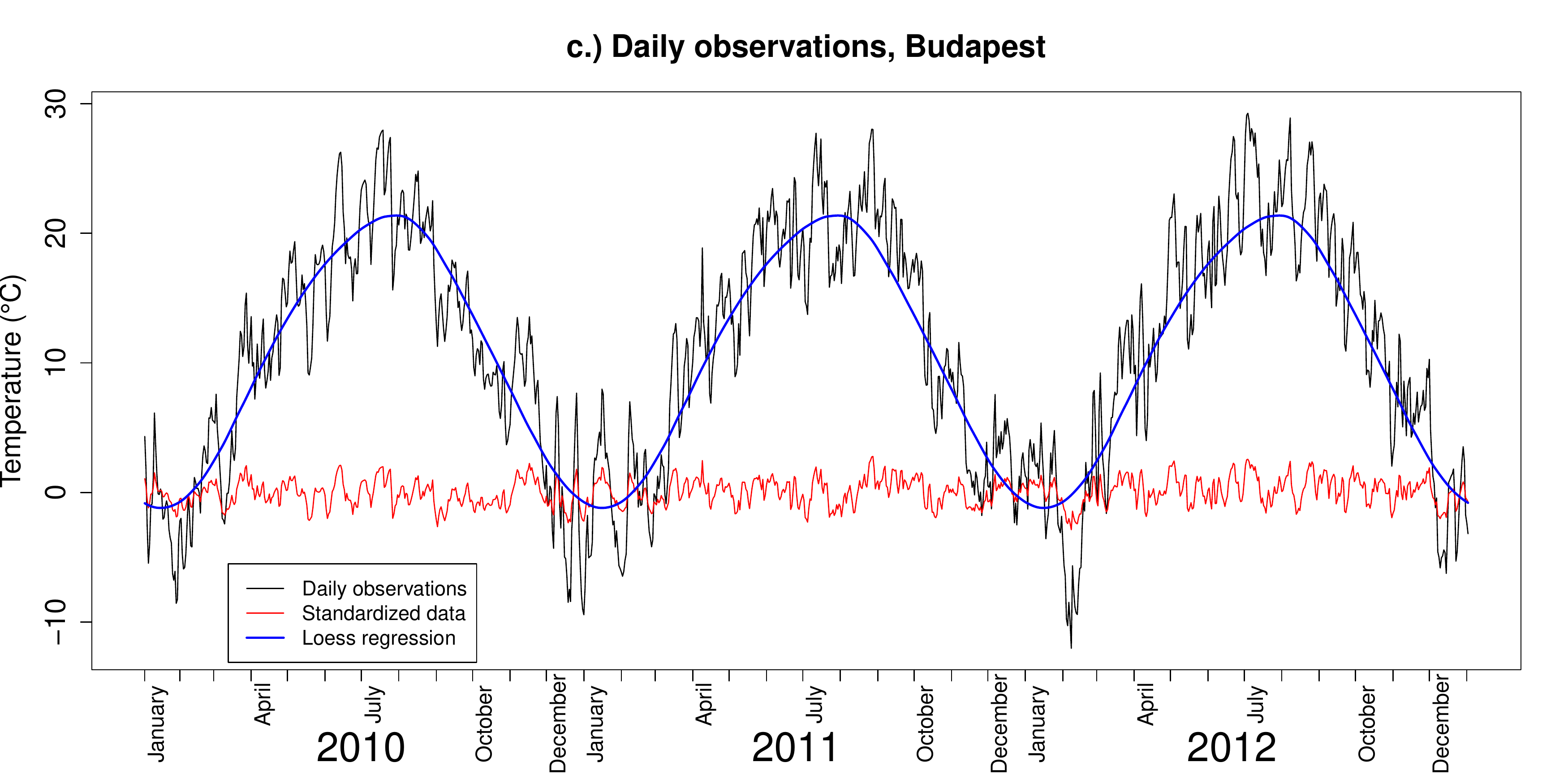} 
	}
	\caption{a.) The daily averages (annual cycle) and the smoothing loess regression; the daily standard deviations (annual cycle) and the smoothing loess regression; c.) the original and the standardized data for grid point Budapest.}
	\label{daily}
\end{figure}

In order to reduce the effect of the outliers of our results, we have finally computed the ten-days averages of the $\tilde{x}$ values. There is a slight but significant upward linear trend in the data, but we did not remove it, as one of our main aims was to detect the changes in the dependencies of the investigated sites -- and these should be based on the original (standardised) deviations, as constructed above.

In the next step, we examined the fixed grid point Budapest paired with other grid points of the database. Using the Akaike information criterion, we chose the lags of the most appropriate among vector autoregressions to model our data pairs. Despite the adjusted $R^2$ values being rather low (around 10\%), the Ljung-Box Q-test, the Breusch-Godfrey test and the test of \citet{ky} could not detect the presence of further serial dependence that has not been included in the VAR model. Table \ref{port} contains these results.

\begin{table}[!h]
	\caption{Results of tests checking for serial dependence between the estimated residuals (p-values) and the chosen lag of the VAR model by Akaike criterion.}
	\label{port}
	\small
\begin{tabular}{|l|c||c|c|c|}
	\hline
	\multirow{2}{*}{Pairs of selected grid points} & Optimal & Ljung- & Breusch- & Genest-Rémillard- \\
	& lag & Box & Godfrey & Kojadinovic-Yan  \\
   \hline
	Budapest \& Sopron & 3 & 0.271 & 0.324 & 0.259 \\
	Budapest \& Apatovac & 1 & 0.080 & 0.463 & 0.120 \\
	Budapest \& Zaránd Mountains & 9 & 0.174 & 0.017 & 0.149 \\
	Budapest \& Nyíregyháza & 4 & 0.174 & 0.474 &  0.258 \\
	Budapest \& Püspökhatvan & 4 & 0.155 & 0.539 & 0.276 \\
	\hline
\end{tabular}     
\end{table}
     
Our main goal is to detect if there is a significant change in the dependence structure of the data. We separated the pairs of points into two parts -- the first part corresponds to the first 31.5 years' observations and the second part to the second 31.5 years' observations. For five selected pairs of grid points, we wanted to test the null hypothesis if the copula of first half of the sample is equal to the copula of the second half of the sample.

Table \ref{obl} and Figure \ref{blszam} depict the optimal block lengths obtained from solving equation (\ref{bopt2}).
The second column of Table \ref{obl} and the red line of Figure \ref{bopt2} show the trace of the covariance matrices of the mean, calculated from the fitted VAR(1) models, multiplicated by 1186 -- the half of the original sample size. We can see on the left figure that the trace of the covariance matrix of the mean is not monotone in the neighborhood of the optimum, so we have to be cautious. The trace function (black line in Figure \ref{bopt2}) is always continuous, but not necessarily differentiable, resulting from the construction of our generalized block bootstrap method. 
Generally, we noticed, that as the block sizes tend to be smaller, the trace function is nearer to be  monotonic.
This phenomenon can be explained by the expansion of $\text{Cov}_* (\overline{X}_b^*)$ in formula (\ref{blvar}): if the block size is relatively small compared to the sample size, then the first and second terms are much more dominant over the third part -- which contains the effect of the remainder block size. We got pretty small, 6.51 optimal block size for the pairs Budapest \& Zaránd Mountains and 20.43 for Budapest \& Sopron. In case of the five selected pairs, as many as 7 iterations were always enough to solve  equation (\ref{bopt2}). We have to mention that there exist some pairs of grid points -- especially at the southern part of the Carpathian Basin --, for which equation (\ref{bopt2}) is not solvable.

\begin{table}[!h]
\caption{Optimal block length for the first half of the samples for the five selected pairs of grid points.}
\label{obl}
\small
\begin{center}
\begin{tabular}{|l|c|c|c|}
	\hline
	\multirow{2}{*}{Pairs of selected grid points} & \multirow{2}{*}{$n\cdot \text{tr}\left( \text{Cov} ( \overline{X}_{\varr} ) \right)$} & Optimal & Number of \\
	& & block size  & iterations \\
	\hline
	Budapest \& Sopron &  1.87 & 20.43 & 4 \\
	Budapest \& Apatovac &  1.92 &   11.23 & 3  \\
	Budapest \& Zaránd Mountains & 2.00 &  {\color{white} 1}6.51 & 3 \\
	Budapest \& Nyíregyháza & 1.84 &  {\color{white} 1}9.34 & 3 \\
	Budapest \& Püspökhatvan &  1.92 & 32.35 & 7 \\
	\hline
\end{tabular}  
\end{center}
\end{table}

\begin{figure}[!h]
	\centering{
		    \includegraphics[width=6cm]{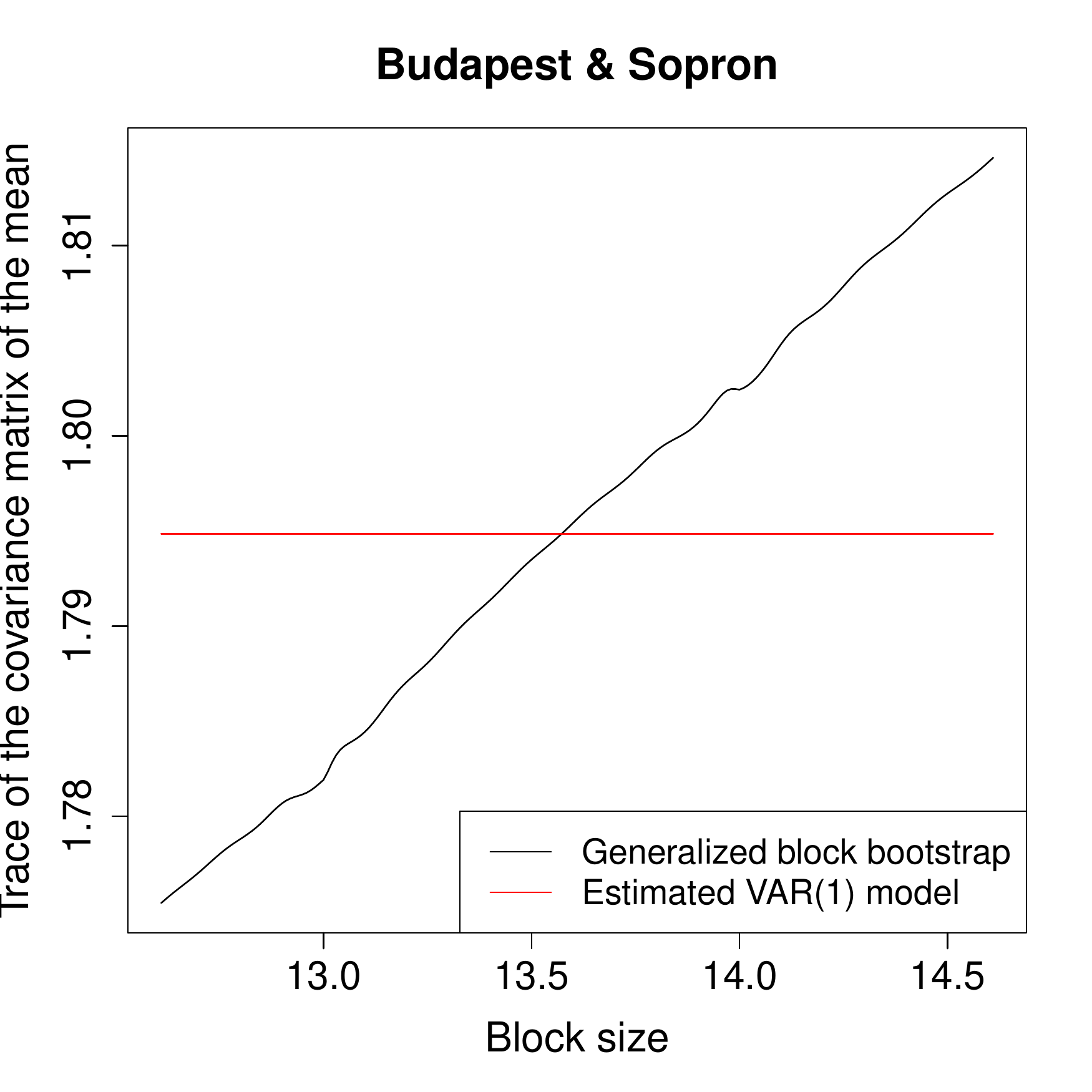}  \includegraphics[width=6cm]{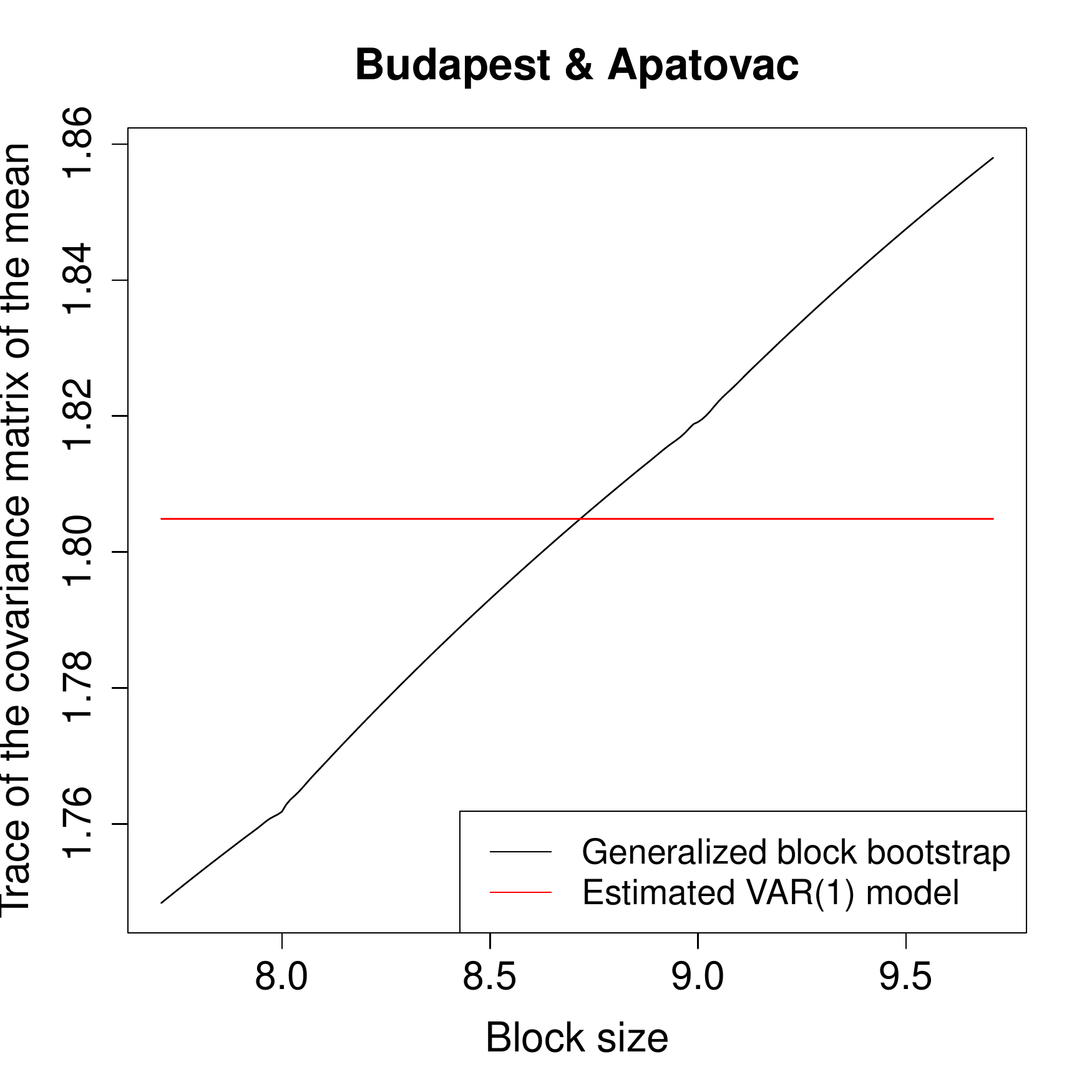} 
	}
	\caption{The trace of the covariance matrix of the mean for two selected grid points (first half of the sample).}
	\label{blszam}
\end{figure}

The last step was conducting the copula homogeneity test described in Section \ref{ChCop}. Using the optimal block size, we generated  boostrap samples via the generalized block bootstrap method for the first half of the sample. With the empirical copulas of these boostrap samples and the empirical copula of the second part of the original sample, we can calculate the test statistic $S_{n,m}$. In our case $n=m=1186$. In order to get accurate $p$-values, we used 2000 repetitions. Table \ref{cophyp} contains the results of the homogeneity test. The dependence structure proved to be different in the first 31.5 years at the two pairs Budapest \& Apatovac and Budapest \& Zaránd Mountains. Figure \ref{EmpCop} depicts the standarized observations and their copula of the pair Budapest \& Apatovac. We can see that preudo-observations are apparently somewhat different, and the test also detected deviance between the two copulas.
 
\begin{table}[!h]
	\caption{$p$-values of the copula homogeneity test, based on 2000 simulations.}
	\label{cophyp}
	\small
	\begin{center}
		\begin{tabular}{|l|c|}
			\hline
			\multirow{2}{*}{Pairs of selected grid points} & \multirow{2}{*}{$p$-values}  \\
			& \\
			\hline
			Budapest \& Sopron &  0.064  \\
			Budapest \& Apatovac &  0.028  \\
			Budapest \& Zaránd Mountains &  0.034  \\
			Budapest \& Nyíregyháza & 0.116 \\
			Budapest \& Püspökhatvan &  0.848 \\
			\hline
		\end{tabular}  
	\end{center}
\end{table}

\begin{figure}[!h]
	\centering{
		\includegraphics[width=6cm]{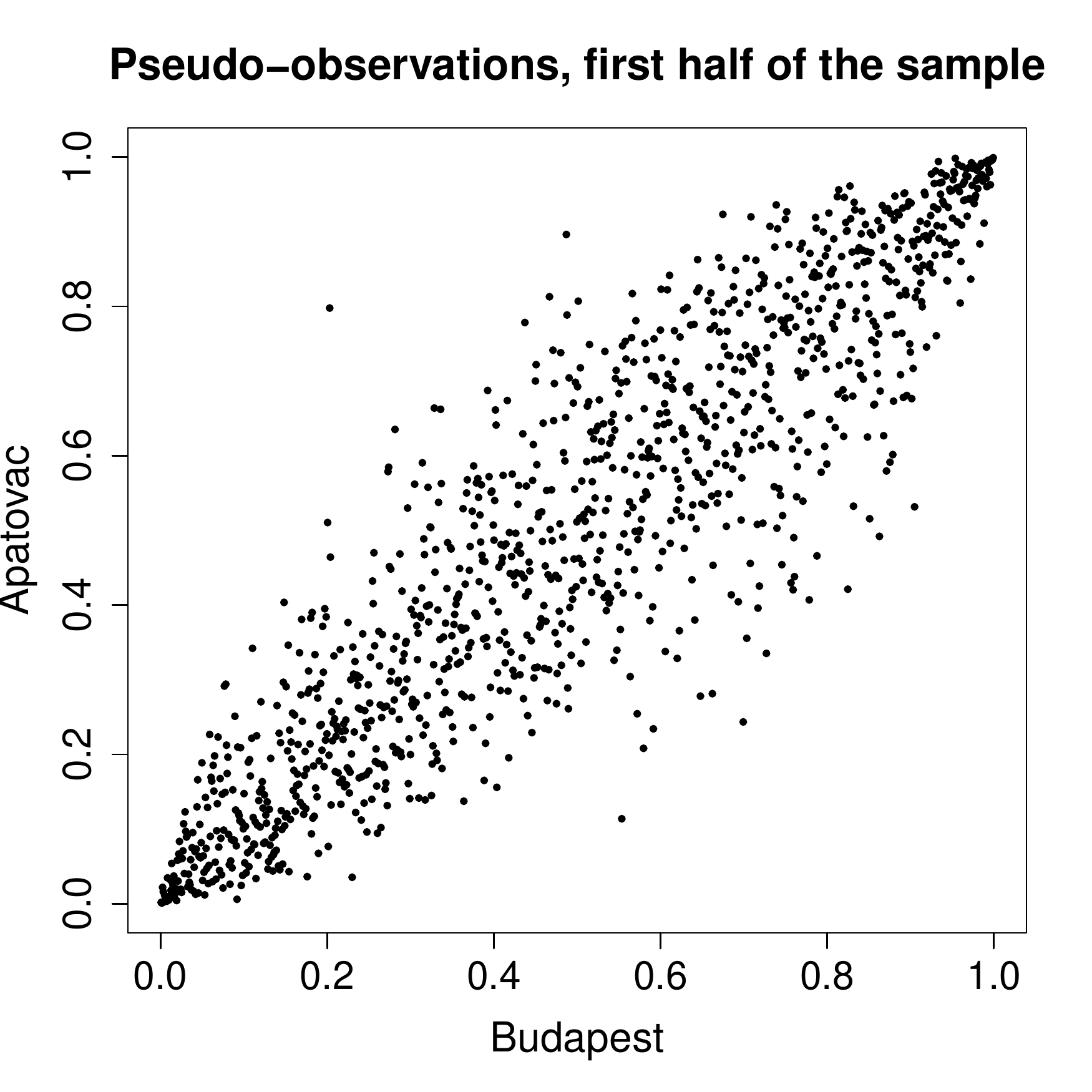}  \includegraphics[width=6cm]{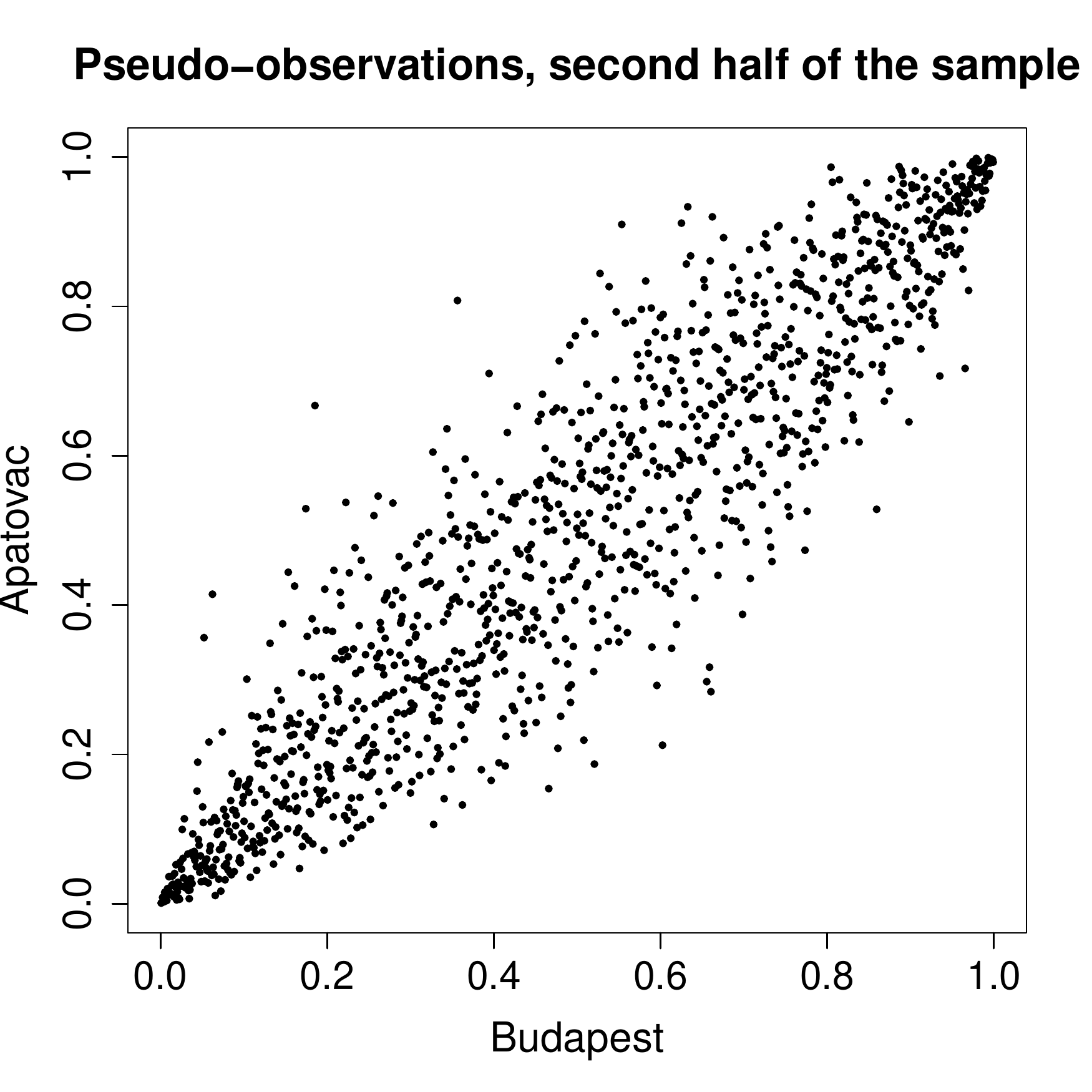} 
	}
	\caption{The pseudo-observations of Budapest \& Apatovac for the first and the second half of the sample.}
	\label{EmpCop}
\end{figure}
    % alkalmazások

%%% 6. fejezet %%%

\section{Conclusions}

We can summarize our results as follows.

There are two major points to be mentioned. First, we proposed a simple generalisation of the block bootstrap methodology, which fits naturally to the existing algorithms, and which helps to overcome the problem of discreteness in the usual block size. 
%Second, since the serial dependence has a substantial impact on the critical values of the used tests, 
%it is important to estimate this effect. The applied block bootstrap methodology was a promising step towards this aim, we were able to estimate the effect of dependence by this tool. 
Second, we have found some significant changes in the dependence structure between the standardised temperature values of pairs of stations within the Carpathian Basin. The direction 
of this change may be worth for investigating, as this may lead to a better understanding of the processes of our climate. 

The proposed generalized block bootstrap method can easily be applied to any other problem, where the block size plays an important role, as all block length determining algorithms give a real number as estimated block size.

It is an interesting question, to which models the proposed block size determining method can be successfully applied. We have checked by simulation that the method is reasonably stable for the VAR models, presented in the paper. For nonlinear time series we might need more observations to fit, which is similarly reliable as the one, presented in our paper.

\vspace{3mm}
\noindent {\large \textbf{Acknowledgements}} \\[1mm]
We acknowledge the E-OBS dataset from the EU-FP6 project ENSEMBLES (\url{http://ensembles-eu.metoffice.com}) and the data providers in the ECA\&D project (\url{http://www.ecad.eu}).

    % summary

\bibliographystyle{plainnat}
\bibliography{bibl_0412} 

\end{document}